\documentclass[oneside,english]{myart}
\usepackage[T1]{fontenc}
\usepackage[latin9]{inputenc}
\usepackage{amssymb}

\makeatletter
\usepackage{natbib}
\usepackage{amssymb}

\makeatother

\usepackage{babel}

\begin{document}
\begin{frontmatter}

\title{Cayley-Dicksonia Revisited}

\author{U. Merkel}

\address{D-70569 Stuttgart, Universitätsstr. 38, merkel.u8@googlemail.com\\
From a 1981 manuscript, unpublished then, now re-typeset for the net
community of the hypercomplex, with my compliments}

\begin{abstract}
In the theory of the hypercomplex, the laws governing the algebra
are based on units that are naturally associated with an orthogonal
vector space, a requirement that is far from mandatory in many algebraic
formulations arising in the context of the reals or the complex numbers.$^{1}$%
\thanks{To name but two prominent examples: algebras with units 1,$\Phi$
in the reals, and 1,$\varpi$ in the complex numbers, where $\Phi=(1+\sqrt{5})/2$
and $\varpi=(1+\sqrt{3}\, i)/2\,,$ respectively.%
}\,\,In this article the complementing view is held, in that the
laws of hypercomplex algebra are recast in terms of quite generally
posited units. Proceeding in this manner, a generalized form of the
Cayley-Dickson process is examined. The representations given are
regular bimodular; the resulting matrices are standard except they
are allowed nonstandard multiplication for noncommutative matrix elements.
\end{abstract}
\end{frontmatter}

\section{Introduction}

Among the infinite set of all algebras, there is one class that stands
out for its members -- the algebra of real numbers, the algebra of
the complex numbers, the algebra of quaternions, and the octonion
(or Cayley-)algebra. These algebras, known as hypercomplex algebras,
share some unique properties which, in part or altogether, are absent
from either other algebras or those venturing beyond.

\emph{Frobenius' theorem}: Each associative division algebra is isomorphic
to either the algebra of real numbers, the algebra of complex numbers
or the quaternion algebra {[}1-3]; for the more general class of alternative
division algebras, the algebra of octonions has to be included.

\emph{Hurwitz' theorem}: Each normed composition algebra with unit
element is isomorphic to one of the algebras listed above {[}4-7]. 

The nonassociativity and/or nonalternativity of an algebra gets in
the way of matrix representation in its usual form. The very uniqueness
of properties of hypercomplex algebras, however, calls for a continuation
of the matrix-representation approach, and in \emph{4.2} we present
matrices with nonstandard multiplication which allow accommodating
nonassociativity or nonalternativity.

As to the applicability of (non-)associative or (non-)alternative
hypercomplex algebras in physics, the reader is referred to {[}8].

\section{Choice of units}

It is common ground in the theory of the hypercomplex to start from
a vector space from which ensue units that are orthogonal to each
other. To free ourselves from the limitations of this approach, here,
we seek to proceed in a different manner.

Hypercomplex algebra shall be defined as a linear algebra over the
field of real numbers $\mathbb{R}$ (or complex numbers $\mathbb{C}$),
with general element (summation assumed) \begin{equation}
\qquad X=x_{0}e_{0}+x_{i}e_{i},\qquad x_{0},x_{i}\in\mathbb{R}\;\mathrm{\left(or\;\mathbb{C}\right)}\end{equation}
where $i=1$ for quadratic algebra, $i=1,2,3$ for quaternion algebra,
and $i=1,2,\dots,7$ for octonion algebra.

The units $e_{i}$ satisfy the multiplication rules\begin{eqnarray}
\qquad\begin{array}{l}
e_{i}^{2}\\
e_{i}e_{0}\\
e_{0}e_{0}\end{array} & \begin{array}{l}
=\\
=\\
=\end{array} & \begin{array}{l}
-qe_{0}-pe_{i},\quad p,q\in\mathbb{R\:\mathrm{\left(or\;\mathbb{C}\right)}}\\
e_{0}e_{i}=e_{i},\\
e_{0},\end{array}\end{eqnarray}
which is all what is needed to describe the class of quadratic algebras,
denoted here $\mathbf{C}\left(p,q\right)$. For the quaternion algebras,
$\mathbf{Q}\left(p,q\right)$, and octonion algebras, $\mathbf{O}\left(p,q\right)$,
one further rule has to be added:\begin{equation}
\quad e_{i}e_{j}=\left(\delta_{ij}D+\epsilon_{ijk}\frac{p}{2}\sqrt{-D}-\frac{p^{2}}{4}\right)e_{0}-\frac{p}{2}e_{i}-\frac{p}{2}e_{j}+\epsilon_{ijk}\sqrt{-D}e_{k}\,,\end{equation}
where $D=\frac{p^{2}}{4}-q$, and $\delta_{ij}$ is the Kronecker
symbol and $\epsilon_{ijk}$ the fully anti\-symmetric Levi-Civita
symbol. In a particular basis, we have $\epsilon_{ijk}=1$, with cyclic
triples $ijk=123$, $145$, $176$, $246$, $257$, $347$, $365$
for octonions, and any of these triples for the quaternion subalgebra.

From Eqs. (2) and (3) we get the commutation and anticommutation relations\begin{eqnarray}
\qquad\begin{array}{l}
\left[e_{i},e_{j}\right]\\
\left\{ e_{i},e_{j}\right\} \end{array} & \begin{array}{l}
=\\
=\end{array} & \begin{array}{l}
2\epsilon_{ijk}\sqrt{-D}(\frac{p}{2}e_{0}+e_{k})\,,\\
2\left((\delta_{ij}D-\frac{p^{2}}{4})e_{0}-\frac{p}{2}(e_{i}+e_{j})\right)\end{array}\end{eqnarray}
and the association relation\begin{eqnarray*}
\begin{array}{l}
\left(e_{i,},e_{j},e_{k}\right)\\
\\\\\\\\\end{array} & \begin{array}{l}
=\\
=\\
\\\\\\\end{array} & \begin{array}{l}
\left(e_{i}e_{j}\right)e_{k}-e_{i}\left(e_{j}e_{k}\right)\\
-D\left((-\frac{p}{2}(\delta_{ij}-\delta_{jk}-\epsilon_{ijp}\epsilon_{pkr}+\epsilon_{jkq}\epsilon_{iqr})+\right.\\
\sqrt{-D}(\delta_{iq}\epsilon_{jkq}-\delta_{pk}\epsilon_{ijp}))e_{0}+\\
\delta_{jk}e_{i}-\delta_{ij}e_{k}+\\
\left.(\epsilon_{ijp}\epsilon_{pkr}-\epsilon_{jkq}\epsilon_{iqr})e_{r}\right)\end{array}\end{eqnarray*}
(which, since $\delta_{iq}\epsilon_{jkq}-\delta_{pk}\epsilon_{ijp}$
vanishes identically)\begin{eqnarray}
\qquad\qquad\begin{array}{l}
\\\\\\\end{array} & \begin{array}{l}
=\\
\\\\\end{array} & \begin{array}{l}
-D\left(-\frac{p}{2}(\delta_{ij}-\delta_{jk}-\epsilon_{ijp}\epsilon_{pkr}+\epsilon_{jkq}\epsilon_{iqr})e_{0}+\right.\\
\delta_{jk}e_{i}-\delta_{ij}e_{k}+\\
\left.(\epsilon_{ijp}\epsilon_{pkr}-\epsilon_{jkq}\epsilon_{iqr})e_{r}\right)\,.\end{array}\end{eqnarray}
The quaternion units $e_{1},e_{2},e_{3}$ are noncommutative, although
associative and alternative; the octonion units $e_{1},e_{2},\dots,e_{7}$
are noncommutative, nonassociative, and alternative.

Due to linearity we can calculate products, commutators and associators
for arbitrary quaternions and octonions of the form (1). For an arbitrary
hypercomplex element $X$, the \emph{conjugate element $\bar{X}$
}is defined by\begin{eqnarray}
\qquad\qquad\begin{array}{l}
\bar{X}\\
\bar{e}_{i}\end{array} & \begin{array}{l}
=\\
=\end{array} & \begin{array}{l}
x_{0}e_{0}+x_{i}\bar{e}_{i}\,,\\
-pe_{0}-e_{i}\,.\end{array}\end{eqnarray}
The conjugation mapping $X\mapsto\bar{X}$ is an \emph{involution},
i.e. $\bar{\bar{X}}=X\,$, $\overline{XY}=\bar{Y}\bar{X}\,$. 

Also, the \emph{norm} of $X$ is \begin{equation}
\mathrm{\qquad\qquad N}(X)=X\bar{X}=\bar{X}X=C_{ij}x_{i}x_{j}e_{0}\,,\end{equation}
with symmetric quantities\[
\qquad\begin{array}{l}
C_{00}=1,\: C_{11}=C_{22}=\dots=C_{77}=q\,,\\
C_{01}=C_{10}=C_{02}=C_{20}=\dots=C_{07}=C_{70}=-\frac{p}{2}\,,\\
C_{12}=C_{21}=C_{13}=C_{31}=\dots=C_{67}=C_{76}=\frac{p^{2}}{4}\,,\end{array}\]
and the inverse element $X_{}^{-1}$ is\begin{equation}
\qquad X^{-1}=\bar{X}/\mathrm{N}(X),\quad X^{-1}X=XX^{-1}=1\cdot e_{0}\,.\end{equation}

In order that a quadratic form (7) becomes positive definite, the
well known Hurwitz criterion requires that the main minors of the
associated symmetric matrix be all positive and real. In the case
of the matrix $C$, these are simply non-negative powers of $-D=q-\frac{p^{2}}{4}\,$.
Hence we can say that algebras $\mathbf{C}(p,q)$, $\mathbf{Q}(p,q)$
and $\mathbf{O}(p,q)$ are \emph{division algebras} for $p,q$ real
and $-D$ positive, with nonisomorphic companions when $p$ and $q$
render the norm (7) degenerate or non-definite. To be precise, $X=\bar{B}A/\mathrm{N}(B)$
is the solution of the equation $BX=A$, and $Y=A\bar{B}/\mathrm{N}(B)$
is the solution to $YB=A$. It is easily verified that the norm (7)
satisfies the decomposition property $\mathrm{N}(XY)=\mathrm{N}(X)\mathrm{N}(Y)$.
Hence hypercomplex algebras, with dimensions necessarily specified
as 1,2,4 or 8, are composition algebras. When $\mathrm{N}(XY)$ is
nondegenerate and positive definite, the norm serves as a definition
of the scalar product $(X,Y)$, and the algebra becomes a normed algebra.

\section{The Cayley-Dickson process}

Let $A$ be a $2^{n}$-dimensional linear algebra over $F$ ($\mathbb{R}$
or $\mathbb{C}$) with an identity element $e_{0}$ and involution
$a\mapsto\bar{a\,}$:\[
\quad\begin{array}{l}
\forall a\in A,\:\exists\,\bar{a}\in A\,:\; a+\bar{a},\:\bar{a}a\left(=a\bar{a}\right)\in F,\:\bar{\bar{a}}=a\,;\\
\forall a,b\in A\,:\;\overline{ab}=\bar{b}\bar{a}\,.\end{array}\]
Then by the Cayley-Dickson process {[}14] a $2^{n+1}$-dimensional
linear algebra $B$ can be formed with identity element and involution
and $A\subset B$ as subalgebra:

$\vartriangleright$\,\,\,\,\,\,\,\,\,\,%
\begin{minipage}[t][1\totalheight]{0.925\columnwidth}%
for all possible ordered pairs $\left(a_{1},a_{2}\right)$ formed
of $a_{1},a_{2}\in A$, addition and scalar multiplication is defined
componentwise and multiplication by the formula%
\end{minipage}

\begin{eqnarray}
\;\begin{array}{l}
\left(a_{1,},a_{2}\right)\left(a_{3},a_{4}\right)\\
\\\\\end{array} & \begin{array}{l}
\!\!=\!\!\\
\!\!\!\!\\
\\\end{array} & \begin{array}{l}
(a_{1}a_{3}-\frac{p}{2}\left[a_{1},a_{4}\right]-\frac{p}{2}a_{2}\left(a_{3}-\bar{a}_{3}\right)-q\,\bar{a}_{4}a_{2}+\frac{p^{2}}{2}\left[a_{2},a_{4}\right],\!\!\!\!\!\\
\quad a_{4}a_{1}+a_{2}\bar{a}_{3}-\frac{p}{2}\left(a_{2}\bar{a}_{4}+a_{4}a_{2}\right))\\
\:\:\:\: a_{1},a_{2},a_{3},a_{4}\in A\,,\; p,q\in F.\end{array}\end{eqnarray}
\,\,\,\,\,\,\,\,\,\,\,\,\,\,%
\begin{minipage}[t][1\totalheight]{0.925\columnwidth}%
The set of all pairs $\left(a_{1},a_{2}\right)=b$ forms a $2^{n+1}$-dimensional
algebra $B$ over $F$. The identity element of $A$, $e_{0}$, is
also the identity of $B$, and will be denoted $\left(e_{0},0\right)$.
The set of all elements $\left(a,0\right)$ , $a\in A$, forms a subalgebra
of $B$ isomorphic to $A$. Adjoined to $A$ is the special element
$\tilde{e}=\left(0,e_{0}\right)$ by which we can write for an element
of $B\,$:%
\end{minipage}\\
\begin{equation}
\qquad\qquad b=a_{1}+a_{2}\tilde{e\,;}\;\end{equation}

$\qquad\;\,$multiplication is defined by

\begin{eqnarray}
\begin{array}{l}
\,\left(a_{1}\!+a\tilde{e}\right)\left(a_{3}\!+a_{4}\tilde{e}\right)\\
\\\\\end{array} & \begin{array}{l}
=\\
\\\\\end{array} & \begin{array}{l}
a_{1}a_{3}-\frac{p}{2}\left[a_{1},a_{4}\right]-\frac{p}{2}a_{2}\left(a_{3}-\bar{a}_{3}\right)-q\,\bar{a}_{4}a_{2}+\\
\frac{p^{2}}{2}\left[a_{2},a_{4}\right]+(a_{4}a_{1}+a_{2}\bar{a}_{3}-\frac{p}{2}\left(a_{2}\bar{a}_{4}+a_{4}a_{2}\right))\tilde{e}\,,\\
a_{1},a_{2},a_{3},a_{4}\in A\,,\; p,q\in F,\end{array}\end{eqnarray}

$\qquad\;$and the involution mapping $b\mapsto\bar{b}$ is defined
by\begin{equation}
\qquad\qquad\bar{b}=\bar{a}_{1}-\frac{p}{2}\left(a_{2}+\bar{a}_{2}\right)-a_{2}\tilde{e}\,.\end{equation}
\[
\qquad\qquad\qquad\qquad\qquad\qquad\qquad\qquad\qquad\qquad\square\]

Setting out from the real numbers $\mathbb{R}$, the algebras that
result at the first step are $\mathbf{C}(p,q)$ with general element\begin{equation}
\qquad\mathbf{c}=x_{0}e_{0}+x_{1}e_{1}\,,\qquad e_{1}e_{1}=-qe_{0}-pe_{1}\,,\end{equation}
where $x_{0},x_{1}\in\mathbb{R}$, provided that $p,q\in\mathbb{R}$.
The norm is\begin{equation}
\mathrm{\qquad N}(c)=\mathbf{c\bar{c}}=\mathbf{\bar{c}c}=x_{0}^{2}-px_{0}x_{1}+qx_{1}^{2}.\end{equation}
Taking values $p=0,\; q=\pm1,0$, three important nonisomorphic algebras
are seen to emerge from class $\mathbf{C}(p,q)$: 

\begin{enumerate}
\item complex numbers $\mathbf{C}(0,1)\simeq\mathbb{C};$
\item split complex numbers $\mathbf{C}(0,-1);$
\item dual numbers $\mathbf{C}(0,0)$.
\end{enumerate}
The transformations $\mathbf{\mathbf{C}(\mathrm{0,1})\mapsto C}(-2\cos\frac{\pi}{2}k,1)$,
$\mathbf{\mathbf{C}(\mathrm{0,-1})\mapsto C(\mathrm{-1-\mathit{i^{\mathrm{2}\mathit{k}}},\mathit{i}^{2\mathit{k}}})}$
(where in the latter case $i=\sqrt{-1}$, hence $p,q,x'_{0},x'_{1}\in\mathbb{C}$)
yield what may be termed periodic algebras,%
\footnote{see appendix A%
} a special subclass each, each of whose members is well-defined except
for poles at $k=2n$. 

Proceeding from $\mathbf{C}(p,q)$, we have at the second step algebras
$\mathbf{Q}(p,q)$, with general element\begin{eqnarray}
\begin{array}{l}
\mathbf{\:\: q}\\
\\\mathbf{\:\:\bar{q}}\\
\mathrm{N(\mathbf{q})}\\
\\\end{array} & \begin{array}{l}
=\\
\\=\\
=\\
\\\end{array} & \begin{array}{l}
\mathbf{c}_{0}+\mathbf{c}_{2}e_{2}=\alpha_{0}e_{0}+\alpha_{1}e_{1}+\alpha_{2}e_{2}+\alpha_{3}e_{3},\;\mathbf{c}_{0},\mathbf{c}_{2}\in\mathrm{\mathbf{C}}(p,q);\\
\mathbf{c}_{0}=x_{0}e_{0}+x_{1}e_{1},\;\mathbf{c}_{2}=x_{2}e_{0}+x_{3}e_{1},\; x_{0},x_{1},x_{2},x_{3}\in\mathbb{R};\\
\mathbf{\bar{\mathbf{c}}}_{0}-\frac{p}{2}(\mathbf{c}_{2}+\bar{\mathbf{c}}_{2})-\mathbf{c}_{2}e_{2}=\alpha_{0}e_{0}+\alpha_{1}\bar{e}_{1}+\alpha_{2}\bar{e}_{2}+\alpha_{3}\bar{e}_{3};\\
\bar{\mathbf{q}}\mathbf{q}=\mathbf{q}\bar{\mathbf{q}}=\alpha_{0}^{2}-p\alpha_{0}(\alpha_{1}+\alpha_{2}+\alpha_{3})\:+\\
\frac{p^{2}}{2}(\alpha_{1}\alpha_{2}+\alpha_{1}\alpha_{3}+\alpha_{2}\alpha_{3})+q(\alpha_{1}^{2}+\alpha_{2}^{2}+\alpha_{3}^{2})\,;\end{array}\end{eqnarray}
whence we validate the relations\begin{eqnarray}
\qquad\begin{array}{l}
\alpha_{0}\\
\alpha_{1}\\
\alpha_{3}\end{array} & \begin{array}{l}
=\\
=\\
=\end{array} & \begin{array}{l}
x_{0}+\frac{p}{2}(\sqrt{-D}-\frac{p}{2})x_{3},\\
x_{1}-\frac{p}{2}x_{3},\qquad\alpha_{2}=x_{2}-\frac{p}{2}x_{3,},\\
\sqrt{-D}x_{3},\end{array}\end{eqnarray}
and also, $e_{2}^{2}=-qe_{0}-pe_{2}.$

The three nonisomorphic algebras that arise at this stage are%
\footnote{were it not for the evidence produced by regular bimodular representation
(the subject of the next section), the existence of simple dual numbers
could be revealed only from the vantage point of quaternions,\emph{
vi}z. $\mathbf{C}(0,0)$ $\subset$ $\mathbf{Q}(0,1)$ $\cup$ $\mathbf{Q}(0,-1)$%
} 

\begin{enumerate}
\item Hamilton's quaternions $\mathbf{Q}(0,1)$;
\item $\frac{1}{4}$-quaternions $\mathbf{Q}(0,-1)$;
\item duodual numbers $\mathbf{Q}(0,0)$.
\end{enumerate}
At step three, proceeding from $\mathbf{Q}(p,q)$, we are led to the
algebras $\mathbf{O}(p,q)$, with general element\begin{eqnarray}
\quad\begin{array}{l}
\mathbf{\:\: o}\\
\\\\\:\:\bar{\mathbf{o}}\\
N(\mathbf{o})\\
\\\end{array} & \begin{array}{l}
=\\
\\\\=\\
=\\
\\\end{array} & \begin{array}{l}
\mathbf{q}_{0}+\mathbf{q}_{4}e_{4}=\beta_{0}e_{0}+\beta_{1}e_{1}+\dots+\beta_{7}e_{7},\\
\mathbf{q}_{0},\mathbf{q}_{4}\in\mathbf{Q}(p,q)\,;\mathbf{q}_{0}=\alpha_{0}e_{0}+\alpha_{1}e_{1}+\dots+\alpha_{3}e_{3},\;\\
\mathbf{q}_{4}=\alpha_{4}e_{0}+\alpha_{5}e_{1}+\dots+\alpha_{7}e_{3};\\
\bar{\mathbf{q}}_{0}-\frac{p}{2}(\mathbf{q}_{4}+\bar{\mathbf{q}}_{4})-\mathbf{q}_{4}e_{4}=\beta_{0}e_{0}+\beta_{1}\bar{e}_{1}+\dots+\beta_{7}\bar{e}_{7};\\
\bar{\mathbf{o}}\mathbf{o}=\mathbf{o}\bar{\mathbf{o}}=\beta_{0}^{2}-p\beta_{0}(\beta_{1}+\dots+\beta_{7})\,+\\
\frac{p^{2}}{2}(\beta_{1}\beta_{2}+\beta_{1}\beta_{3}+\dots+\beta_{6}\beta_{7})+q(\beta_{1}^{2}+\beta_{2}^{2}+\dots+\beta_{7}^{2})\,;\end{array}\end{eqnarray}
Eq. (17) is accompanied by the relations\begin{eqnarray}
\quad\begin{array}{l}
\beta_{0}\\
\beta_{1}\\
\beta_{4}\\
\beta_{5}\\
e_{4}^{2}\end{array} & \begin{array}{l}
=\\
=\\
=\\
=\\
=\end{array} & \begin{array}{l}
\alpha_{0}+\frac{p}{2}(\sqrt{-D}-\frac{p}{2})(\alpha_{5}+\alpha_{6}+\alpha_{7}),\\
\alpha_{1}-\frac{p}{2}\alpha_{5},\;\beta_{2}=\alpha_{2}-\frac{p}{2}\alpha_{6},\;\beta_{3}=\alpha_{3}-\frac{p}{2}\alpha_{7},\\
\alpha_{4}-\frac{p}{2}(\alpha_{5}+\alpha_{6}+\alpha_{7}),\\
\sqrt{-D}\alpha_{5},\;\beta_{6}=\sqrt{-D}\alpha_{6},\;\beta_{7}=\sqrt{-D}\alpha_{7},\\
-qe_{0}-pe_{4}\,.\end{array}\end{eqnarray}

Again, three important nonisomorphic algebras can be discerned: 

\begin{enumerate}
\item Cayley's octonions $\mathbf{O}(0,1)$;
\item $\frac{1}{8}$-octonions $\mathbf{O}(0,-1)$;
\item tridual numbers $\mathbf{\mathbf{O}}(0,0)$.
\end{enumerate}
At each step, up to equivalence within the class involved, we have
one division algebra $(q=1)$, one split algebra $(q=-1)$, and one
algebra with nilideal $(q=0)$, where the problem of isomorphism can
be tackled in the following way: According to Jacobson {[}9], two
Cayley-Dickson algebras $A$ and $A'$ are isomorphic if their bilinear
forms $(X,Y)=\frac{1}{2}\left[\mathrm{N}(X+Y)-\mathrm{N}(X)-\mathrm{N}(Y)\right]$
are equivalent, i.e., if there is a linear mapping $X\mapsto XH$,
$Y\mapsto YH$ of $A$ into $A'$ such that $\mathrm{N}(X)=\mathrm{N'}(XH)$,
$\mathrm{N}(Y)=\mathrm{N'}(YH)$, for all $X,Y\in A$.

The Cayley-Dickson-type normed composition algebras allowed by Frobenius'
and Hurwitz' theorems are exhausted after the first three steps. After
step 3 is also brought to a halt alternativity, since the Cayley-Dickson
process yields alternative algebras only if the initial algebra is
associative; meta-Hurwitz algebras therefore fulfil but much weaker
identities such as the flexibility law $X(YX)=(XY)X$ {[}12-14].

\section{Representation of hypercomplex algebra}

To keep this section brief, we just give three specific examples: 

\begin{enumerate}
\item $\mathbf{Q}(p,q)$-rep by matrices over the field $\mathbb{C}\simeq\mathbf{C}(0,1)$; 
\item $\mathbf{O}(0,1)$-rep by matrices over the quaternionic skew field
$\mathbb{Q\simeq\mathbf{Q}}(0,1)$, and
\item $\mathbf{S}(0,1)$-rep by matrices over the octonionic \emph{aq}-field
(or alternative quasi-field) $\mathbf{\mathbb{O}\simeq}$ $\mathbf{O}(0,1)$.
\end{enumerate}
For examples (2) and (3), rep matrices are nonstandard, although in
the first example they are in perfect compliance with representation
lore {[}10].

\subsection{Representation of ${\mathbf Q}(p,q)$}

Eq. (2) can be regarded as a quadratic equation which resolves into\begin{eqnarray}
\;\begin{array}{l}
(e_{1})_{1,2}\\
\\\end{array} & \begin{array}{l}
=\\
\\\end{array} & \begin{array}{l}
-\frac{p}{2}\pm\sqrt{D},\qquad(D=\frac{p^{2}}{4}-q)\\
p,q\in\mathbb{C}.\end{array}\end{eqnarray}
Thus we have\begin{equation}
\!\!\!\!\!\!\!\!\!\!\!\!\!\begin{array}{l}
(\alpha_{0}e_{0}+\alpha_{1}e_{1}+\alpha_{2}e_{2}+\alpha_{3}e_{3})\cdot e_{1}\simeq(x_{0}e_{0}+x_{1}e_{1},x_{2}e_{0}+x_{3}e_{1})\cdot(e_{1},0)\\
=((x_{0}e_{0}+x_{1}e_{1})e_{1}-\frac{p}{2}(x_{2}e_{0}+x_{3}e_{1})(p+2e_{1})\,,(x_{2}e_{0}+x_{3}e_{1})(-p-e_{1}))\\
\simeq((x_{0}e_{0}+x_{1}e_{1})(-\frac{p}{2}\pm\sqrt{D})+(x_{2}e_{0}+x_{3}e_{1})(\mp p\sqrt{D}),\\
\qquad(x_{2}e_{0}+x_{3}e_{1})(-\frac{p}{2}\mp\sqrt{D}))\\
=(x_{0}e_{0}+x_{1}e_{1},x_{2}e_{0}+x_{3}e_{1})\cdot\left(\begin{array}{cc}
-\frac{p}{2}\pm\sqrt{D} & 0\\
\mp p\sqrt{D} & -\frac{p}{2}\mp\sqrt{D}\end{array}\right);\end{array}\end{equation}
\begin{equation}
\!\!\!\!\!\!\!\!\!\!\!\!\!\begin{array}{l}
(\alpha_{0}e_{0}+\alpha_{1}e_{1}+\alpha_{2}e_{2}+\alpha_{3}e_{3})\cdot e_{2}\simeq(x_{0}e_{0}+x_{1}e_{1},x_{2}e_{0}+x_{3}e_{1})\cdot(0,e_{0})\\
=(-qe_{0}(x_{2}e_{0}+x_{3}e_{1})\,,e_{0}(x_{0}e_{0}+x_{1}e_{1})-p(x_{2}e_{0}+x_{3}e_{1})))\\
\simeq((x_{2}e_{0}+x_{3}e_{1})(-q),x_{0}e_{0}+x_{1}e_{1}+(x_{2}e_{0}+x_{3}e_{1})(-p))\\
\\=(x_{0}e_{0}+x_{1}e_{1},x_{2}e_{0}+x_{3}e_{1})\cdot\left(\begin{array}{cc}
0 & 1\\
-q & -p\end{array}\right);\end{array}\end{equation}
\begin{equation}
\!\!\!\!\!\!\!\!\!\!\!\begin{array}{l}
(\alpha_{0}e_{0}+\alpha_{1}e_{1}+\alpha_{2}e_{2}+\alpha_{3}e_{3})\cdot e_{3}\simeq(x_{0}e_{0}+x_{1}e_{1},x_{2}e_{0}+x_{3}e_{1})\\
\cdot\left(\frac{p^{2}-2p\sqrt{-D}}{4\sqrt{-D}}e_{0}+\frac{p}{2\sqrt{-D}}e_{1}\,,\frac{p}{2\sqrt{-D}}e_{0}+\frac{1}{\sqrt{-D}}e_{1}\right)\\
=((x_{0}e_{0}+x_{1}e_{1})(\frac{p^{2}-2p\sqrt{-D}}{4\sqrt{-D}}e_{0}+\frac{p}{2\sqrt{-D}}e_{1})+q(\frac{p}{2\sqrt{-D}}e_{0}+\frac{1}{\sqrt{-D}}e_{1})\\
\cdot(x_{2}e_{0}+x_{3}e_{1})-\frac{p}{2}(x_{2}e_{0}+x_{3}e_{1})(\frac{p^{2}}{2\sqrt{-D}}e_{0}+\frac{p}{\sqrt{-D}}e_{1})\,,\\
(\frac{p}{2\sqrt{-D}}e_{0}+\frac{1}{\sqrt{-D}}e_{1})(x_{0}e_{0}+x_{1}e_{1})+(x_{2}e_{0}+x_{3}e_{1}))\\
\simeq((x_{0}e_{0}+x_{1}e_{1})(-1\mp i)\frac{p}{2}+(x_{2}e_{0}+x_{3}e_{1})(\pm i)(\frac{p^{2}}{2}-q),\\
\qquad(x_{0}e_{0}+x_{1}e_{1})(\mp i)+(x_{2}e_{0}+x_{3}e_{1})(-1\pm i)\frac{p}{2})\\
=(x_{0}e_{0}+x_{1}e_{1},x_{2}e_{0}+x_{3}e_{1})\cdot\left(\begin{array}{cc}
(-1\mp i)\frac{p}{2} & \mp i\\
\pm i(\frac{p^{2}}{2}-q) & (-1\pm i)\frac{p}{2}\end{array}\right).\\
\mathrm{Of\: course,\:\:}\: e_{0}\leftrightarrow\left(\begin{array}{cc}
1 & 0\\
0 & 1\end{array}\right).\end{array}\end{equation}
\subsection{Representation of ${\mathbf O}(0,1)$ and ${\mathbf S}(0,1)$} 

Proceeding further, we let nonstandard multiplier positions assume
the multi\-plication of rep matrices\begin{equation}
\begin{array}{cc}
\left(\begin{array}{cc}
a_{11} & a_{12}\\
a_{21} & a_{22}\end{array}\right)\cdot\left(\begin{array}{cc}
b_{11} & b_{12}\\
b_{21} & b_{22}\end{array}\right) & \!\!\!\!\!\!\!\!\!\!\!\!\!\!\!\!\!\!\!\!\!\!\!\!\!\!\!\!\!\!\!\!\!\!\!\!\!\!\!\!\!\!\!\!\!\!\!\!\!\!=\left(\begin{array}{cc}
c_{11} & c_{12}\\
c_{21} & c_{22}\end{array}\right)\\
 & :=\left(\begin{array}{cc}
a_{11}b_{11}+b_{21}a_{12} & \quad b_{12}a_{11}+a_{12}b_{22}\\
b_{11}a_{21}+a_{22}b_{21} & \quad a_{21}b_{12}+b_{22}a_{22}\end{array}\right)\end{array}\end{equation}
(a deviation entailing failure of associativity for noncommutative
entries, as well as failure of alternativity for nonassociative entries).
Writing $1,i,j,k$ for the units of the quaternionic skew field $\mathbb{Q}$
defined by\begin{equation}
i^{2}=j^{2}=k^{2}=ijk=-1,\end{equation}
we observe\begin{equation}
\!\begin{array}{l}
(\beta_{0}e_{0}+\beta_{1}e_{1}+\beta_{2}e_{2}+\dots+\beta_{7}e_{7})\cdot e_{1}\\
\simeq(\alpha_{0}e_{0}+\alpha_{1}e_{1}+\alpha_{2}e_{2}+\alpha_{3}e_{3},\:\alpha_{4}e_{0}+\alpha_{5}e_{1}+\alpha_{6}e_{2}+\alpha_{7}e_{3})\cdot(e_{1},0)\\
=((\alpha_{0}e_{0}+\alpha_{1}e_{1}+\alpha_{2}e_{2}+\alpha_{3}e_{3})e_{1},(\alpha_{4}e_{0}+\alpha_{5}e_{1}+\alpha_{6}e_{2}+\alpha_{7}e_{3})(-e_{1}))\\
\simeq((\alpha_{0}e_{0}+\alpha_{1}e_{1}+\alpha_{2}e_{2}+\alpha_{3}e_{3})i,(\alpha_{4}e_{0}+\alpha_{5}e_{1}+\alpha_{6}e_{2}+\alpha_{7}e_{3})(-i))\\
=(\alpha_{0}e_{0}+\alpha_{1}e_{1}+\alpha_{2}e_{2}+\alpha_{3}e_{3},\:\alpha_{4}e_{0}+\alpha_{5}e_{1}+\alpha_{6}e_{2}+\alpha_{7}e_{3})\cdot\left(\begin{array}{cc}
i & 0\\
0 & -i\end{array}\right);\end{array}\end{equation}
\begin{equation}
\!\begin{array}{l}
(\beta_{0}e_{0}+\beta_{1}e_{1}+\beta_{2}e_{2}+\dots+\beta_{7}e_{7})\cdot e_{2}\\
\simeq(\alpha_{0}e_{0}+\alpha_{1}e_{1}+\alpha_{2}e_{2}+\alpha_{3}e_{3},\:\alpha_{4}e_{0}+\alpha_{5}e_{1}+\alpha_{6}e_{2}+\alpha_{7}e_{3})\cdot(e_{2},0)\\
=((\alpha_{0}e_{0}+\alpha_{1}e_{1}+\alpha_{2}e_{2}+\alpha_{3}e_{3})e_{2},(\alpha_{4}e_{0}+\alpha_{5}e_{1}+\alpha_{6}e_{2}+\alpha_{7}e_{3})(-e_{2}))\\
\simeq((\alpha_{0}e_{0}+\alpha_{1}e_{1}+\alpha_{2}e_{2}+\alpha_{3}e_{3})j,(\alpha_{4}e_{0}+\alpha_{5}e_{1}+\alpha_{6}e_{2}+\alpha_{7}e_{3})(-j))\\
=(\alpha_{0}e_{0}+\alpha_{1}e_{1}+\alpha_{2}e_{2}+\alpha_{3}e_{3},\:\alpha_{4}e_{0}+\alpha_{5}e_{1}+\alpha_{6}e_{2}+\alpha_{7}e_{3})\cdot\left(\begin{array}{cc}
j & 0\\
0 & -j\end{array}\right);\end{array}\end{equation}
and so on, the last in the line being\begin{equation}
\!\begin{array}{l}
(\beta_{0}e_{0}+\beta_{1}e_{1}+\beta_{2}e_{2}+\dots+\beta_{7}e_{7})\cdot e_{7}\\
\simeq(\alpha_{0}e_{0}+\alpha_{1}e_{1}+\alpha_{2}e_{2}+\alpha_{3}e_{3},\:\alpha_{4}e_{0}+\alpha_{5}e_{1}+\alpha_{6}e_{2}+\alpha_{7}e_{3})\cdot(0,e_{3})\\
=(e_{3}(\alpha_{4}e_{0}+\alpha_{5}e_{1}+\alpha_{6}e_{2}+\alpha_{7}e_{3}),e_{3}(\alpha_{0}e_{0}+\alpha_{1}e_{1}+\alpha_{2}e_{2}+\alpha_{3}e_{3}))\\
\simeq(k(\alpha_{4}e_{0}+\alpha_{5}e_{1}+\alpha_{6}e_{2}+\alpha_{7}e_{3}),k(\alpha_{0}e_{0}+\alpha_{1}e_{1}+\alpha_{2}e_{2}+\alpha_{3}e_{3}))\\
=(\alpha_{0}e_{0}+\alpha_{1}e_{1}+\alpha_{2}e_{2}+\alpha_{3}e_{3},\:\alpha_{4}e_{0}+\alpha_{5}e_{1}+\alpha_{6}e_{2}+\alpha_{7}e_{3})\cdot\left(\begin{array}{cc}
0 & k\\
k & 0\end{array}\right);\end{array}\end{equation}
and, of course,$\; e_{0}\leftrightarrow\left(\begin{array}{cc}
1 & 0\\
0 & 1\end{array}\right).$

In the same vein, writing $1,i,j,\dots,o$ for the units%
\footnote{as did Graves and Hamilton in their famous correspondence%
} of the octonionic alternative quasi-field $\mathbb{O}$, where\begin{equation}
\begin{array}{c}
i^{2}=j^{2}=\dots=o^{2}=-1=ijk=ilm\\
\;=ion=jln=jmo=klo=knm,\end{array}\end{equation}
the analogue of the previous result is obtained, with nonstandard
rep matrices characterizing the nonalternative sedenionic algebra
$\mathbf{S}(0,1)$,\begin{equation}
\begin{array}{cccc}
e_{0}\leftrightarrow\left(\begin{array}{cc}
1 & 0\\
0 & 1\end{array}\right),\\
e_{1}\leftrightarrow\left(\begin{array}{cc}
i & 0\\
0 & -i\end{array}\right), & e_{2}\leftrightarrow\left(\begin{array}{cc}
j & 0\\
0 & -j\end{array}\right), & \dots, & e_{7}\leftrightarrow\left(\begin{array}{cc}
o & 0\\
0 & -o\end{array}\right),\\
e_{8}\leftrightarrow\left(\begin{array}{cc}
0 & 1\\
-1 & 0\end{array}\right),\\
e_{9}\leftrightarrow\left(\begin{array}{cc}
0 & i\\
i & 0\end{array}\right), & e_{10}\leftrightarrow\left(\begin{array}{cc}
0 & j\\
j & 0\end{array}\right), & \dots, & e_{16}\leftrightarrow\left(\begin{array}{cc}
0 & o\\
o & 0\end{array}\right),\end{array}\end{equation}
\medskip{}
 which also concludes our journey into the realm of hypercomplex algebra. 

The author is indebted to Jens Köplinger whose enthusiasm for the
subject provided the inducement to rifle through piles of dust-gathering
manuscripts and rescue this specimen from oblivion. 

\appendix\section{}

From Euler's formula\begin{equation}
\mathrm{e}^{iz}=\cos z+i\sin z\end{equation}
it follows\begin{equation}
\mathrm{e}^{i\pi/2}=i.\end{equation}
Exponentiating both sides of (A.2) with $t$, one obtains the power
law of $i$:\begin{equation}
i^{t}=\cos\frac{\pi}{2}t+i\sin\frac{\pi}{2}t.\end{equation}
(A.3) can be recast in terms of nontrivial roots of unity, $\omega$.
Substituting $t=k\theta$, it reads\begin{equation}
i^{k\theta}=\cos\frac{\pi}{2}k\theta+i\sin\frac{\pi}{2}k\theta,\end{equation}
whence for $\omega=i^{k}$\begin{equation}
\omega^{\theta}=\cos\frac{\pi}{2}k\theta-\cot\frac{\pi}{2}k\sin\frac{\pi}{2}k\theta+\omega\csc\frac{\pi}{2}k\sin\frac{\pi}{2}k\theta.\end{equation}
The power law in the guise of (A.5) may then be generalized in such
a way that it lends itself to the description of two disjoint subclasses
of periodic algebras $\mathbf{C}(-2\rho^{k}\cos\frac{\pi}{2}k,\rho^{2k})$: 

$\mathbf{Theorem}:$ \emph{The power law\begin{equation}
e_{1}^{\theta}=\rho^{k\theta}\left[e_{0}\left(\cos\frac{\pi}{2}k\theta-\cot\frac{\pi}{2}k\sin\frac{\pi}{2}k\theta\right)+e_{1}\rho^{-k}\csc\frac{\pi}{2}k\sin\frac{\pi}{2}k\theta\right]\end{equation}
} \emph{describes two disjoint subclasses of periodic algebras}, i.e.,
\emph{for} $\rho=1$,\emph{ the subclass of complex periodic algebras}
$\mathbf{C}(-2\cos\frac{\pi}{2}k,1)$, \emph{and for} $\rho=i$\emph{,
the split-complex }sub\emph{class} $\mathbf{C}(-1-i^{2k},i^{2k})$.\emph{ }

\emph{Proof}: Case $\rho=1$ holds trivially by way of identification
$e_{1}\leftrightarrow\omega$. Case $\rho=i$ follows from a tedious
but straightforward application of (A.6): 

Let\begin{eqnarray*}
\begin{array}{c}
e_{1}^{2\theta}\\
\\\end{array} & \begin{array}{c}
=\\
=\end{array} & \begin{array}{c}
i^{2k\theta}\left[e_{0}\left(a_{1}(2\theta)-a_{2}(2\theta)\right)+e_{1}b(2\theta)\right]\\
i^{2k\theta}\left[e_{0}\left(a_{1}(\theta)-a_{2}(\theta)\right)+e_{1}b(\theta)\right]^{2}\end{array}\end{eqnarray*}
or\[
\begin{array}{ccc}
i^{2k\theta}\left[e_{0}\left(\underbrace{\cos^{2}\frac{\pi}{2}k\theta}\right.\right. & \underbrace{-2\cot\frac{\pi}{2}k\cos\frac{\pi}{2}k\theta\sin\frac{\pi}{2}k\theta} & \underbrace{+\left.\cot^{2}\frac{\pi}{2}k\sin^{2}\frac{\pi}{2}k\theta\right)}\\
\quad\qquad\left(a_{1}^{2}(\theta)\right. & -2a_{1}(\theta)a_{2}(\theta) & +\left.a_{2}^{2}(\theta)\right)\end{array}\]
\[
\qquad\begin{array}{ccc}
+e_{1}\left(\underbrace{2i^{-k}\csc\frac{\pi}{2}k\cos\frac{\pi}{2}k\theta\sin\frac{\pi}{2}k\theta}\right.\\
\left(2a_{1}(\theta)b(\theta)\right.\end{array}\]
\[
\qquad\qquad\begin{array}{cc}
\left.\underbrace{-2i^{-k}\cot\frac{\pi}{2}k\csc\frac{\pi}{2}k\sin^{2}\frac{\pi}{2}k\theta}\right) & +e_{1}^{2}\left.\underbrace{i^{-2k}\csc^{2}\frac{\pi}{2}k\sin^{2}\frac{\pi}{2}k\theta}\right]\\
\left.-2a_{2}(\theta)b(\theta)\right) & b^{2}(\theta)\end{array}.\]
By definition, the relation $e_{1}^{2}=-qe_{0}-pe_{1}$ holds for
$\mathbf{C}(p,q)$. Whence we have \[
\begin{array}{llll}
(i) & a_{1}^{2}(\theta)-2a_{1}(\theta)a_{2}(\theta)+a_{2}^{2}(\theta)-qb^{2}(\theta) & := & a_{1}(2\theta)-a_{2}(2\theta),\\
(ii) & 2a_{1}(\theta)b(\theta)-2a_{2}(\theta)b(\theta)-pb^{2}(\theta) & := & b(2\theta).\end{array}\]
Which, using $2a_{1}(\theta)a_{2}(\theta)=a_{2}(2\theta)$ and $2a_{1}(\theta)b(\theta)=b(2\theta)$
, reduces to \[
\begin{array}{llll}
(i)' & a_{1}^{2}(\theta)+a_{2}^{2}(\theta)-qb^{2}(\theta) & := & a_{1}(2\theta),\\
(ii)' & -2a_{2}(\theta)b(\theta)-pb^{2}(\theta) & := & 0.\end{array}\]
Keeping in mind that $p=-1-i^{2k},q=i^{2k}$ was assumed, for $(i)'$
one verifies\begin{eqnarray*}
\begin{array}{c}
\cos^{2}\frac{\pi}{2}k\theta+\cot^{2}\frac{\pi}{2}k\sin^{2}\frac{\pi}{2}k\theta-\csc^{2}\frac{\pi}{2}k\sin^{2}\frac{\pi}{2}k\theta\\
\\\end{array} & \begin{array}{c}
=\\
=\end{array} & \begin{array}{c}
\cos^{2}\frac{\pi}{2}k\theta-\sin^{2}\frac{\pi}{2}k\theta\\
\cos\pi k\theta=a_{1}(2\theta);\end{array}\end{eqnarray*}
 and for $(ii)'$,\[
\!\!\!\!\!\!\!\!\!\!\begin{array}{l}
-2i^{-k}\cot\frac{\pi}{2}k\csc\frac{\pi}{2}k\sin^{2}\frac{\pi}{2}k\theta+i^{-2k}\csc^{2}\frac{\pi}{2}k\sin^{2}\frac{\pi}{2}k\theta+\csc^{2}\frac{\pi}{2}k\sin^{2}\frac{\pi}{2}k\theta\\
=i^{-2k}\csc^{2}\frac{\pi}{2}k(-2i^{k}\cos\frac{\pi}{2}k\sin^{2}\frac{\pi}{2}k\theta+\sin^{2}\frac{\pi}{2}k\theta+i^{2k}\sin^{2}\frac{\pi}{2}k\theta)\\
=i^{-2k}\csc^{2}\frac{\pi}{2}k(-2\cos^{2}\frac{\pi}{2}k\sin^{2}\frac{\pi}{2}k\theta-2i\cos\frac{\pi}{2}k\sin\frac{\pi}{2}k\sin^{2}\frac{\pi}{2}k\theta+\sin^{2}\frac{\pi}{2}k\theta\\
\qquad+\cos^{2}\frac{\pi}{2}k\sin^{2}\frac{\pi}{2}k\theta-\sin^{2}\frac{\pi}{2}k\sin^{2}\frac{\pi}{2}k\theta+2i\cos\frac{\pi}{2}k\sin\frac{\pi}{2}k\sin^{2}\frac{\pi}{2}k\theta)\\
=i^{-2k}\csc^{2}\frac{\pi}{2}k\cdot0=0.\end{array}\]
\[
\quad\qquad\qquad\qquad\qquad\qquad\qquad\qquad\qquad\qquad\qquad\qquad\qquad\qquad\square\]

The periodic algebras\emph{ }$\mathbf{C}(-2\rho^{k}\cos\frac{\pi}{2}k,\rho^{2k})$
and $\mathbf{Q}(-2\rho^{k}\cos\frac{\pi}{2}k,\rho^{2k})$, $\rho=1\vee i$\emph{,
}are uniquely linked to\emph{ }their respective orthogonal forms,
derived from (A.6) by setting $k=1$:\begin{equation}
e_{n}^{\theta}=e_{0}\cos\frac{\pi}{2}\theta+e_{n}\sin\frac{\pi}{2}\theta,\qquad\alpha_{0}e_{0}+\alpha_{n}e_{n}\in\mathbf{Q}(0,1),\end{equation}
$\qquad$\begin{equation}
\begin{array}{c}
e_{n}^{\theta}=(\cos\frac{\pi}{2}\theta+i\sin\frac{\pi}{2}\theta)(e_{0}\cos\frac{\pi}{2}\theta-ie_{n}\sin\frac{\pi}{2}\theta),\\
\qquad\qquad\qquad\qquad\qquad\qquad\alpha_{0}e_{0}+\alpha_{n}e_{n}\in\mathbf{Q}(0,-1),\end{array}\end{equation}
where $n=1,2,3$, $ie_{0}=e_{0}i$, $ie_{n}=e_{n}i$. 

$\mathbf{Theorem}:$ \emph{The periodic algebras }$\mathbf{C}(-2\rho^{k}\!\cos\frac{\pi}{2}k,\rho^{2k})$
and $\mathbf{Q}(-2\rho^{k}\!\cos\frac{\pi}{2}k,\rho^{2k})$, $\rho=1\vee i$\emph{,
provide a dual description of their respective orthogonalizations
}(A.7)\emph{ and }(A.8)\emph{ by way of their specific representations.}

\emph{Proof: }As $\mathbf{C}(p,q)$ forms a subalgebra of $\mathbf{Q}(p,q)$,
it suffices to verify the assertion for $\mathbf{Q}(-2\rho^{k}\cos\frac{\pi}{2}k,\rho^{2k})$.
As was shown in section \emph{4.1}, with $\mathbf{Q}(p,q)$ is associated
the regular representation \[
\begin{array}{cc}
u_{0}\leftrightarrow\left(\begin{array}{cc}
1 & \quad0\\
0 & \quad1\end{array}\right), & u_{1}\leftrightarrow\left(\begin{array}{cc}
-\frac{p}{2}\pm\sqrt{D} & \quad0\\
\mp p\sqrt{D} & \quad-\frac{p}{2}\mp\sqrt{D}\end{array}\right),\\
u_{2}\leftrightarrow\left(\begin{array}{cc}
0 & \quad1\\
-q & \quad-p\end{array}\right), & u_{3}\leftrightarrow\left(\begin{array}{cc}
(-1\mp i)\frac{p}{2} & \quad\mp i\\
\pm i(\frac{p^{2}}{2}-q) & \quad(-1\pm i)\frac{p}{2}\end{array}\right).\end{array}\]
Settling for the lower sign, $\mathbf{Q}(-2\cos\frac{\pi}{2}k,1)$
has representation\[
\!\!\!\!\!\!\!\!\!\!\!\!\!\begin{array}{cc}
(u_{(1)})_{0}\leftrightarrow\left(\begin{array}{cc}
1 & \:0\\
0 & \:1\end{array}\right), & (u_{(1)})_{1}\leftrightarrow\left(\begin{array}{cc}
\cos\frac{\pi}{2}k+i\sin\frac{\pi}{2}k & \quad0\\
2i\cos\frac{\pi}{2}k\sin\frac{\pi}{2}k & \quad\cos\frac{\pi}{2}k-i\sin\frac{\pi}{2}k\end{array}\right),\\
(u_{(1)})_{2}\leftrightarrow\left(\begin{array}{cc}
0 & \:1\\
-1 & \:2\cos\frac{\pi}{2}k\end{array}\right), & (u_{(1)})_{3}\leftrightarrow\left(\begin{array}{cc}
(1-i)\cos\frac{\pi}{2}k & \quad i\\
-i(\cos^{2}\frac{\pi}{2}k-\sin^{2}\frac{\pi}{2}k) & \:(1+i)\cos\frac{\pi}{2}k\end{array}\right),\end{array}\]
and $\mathbf{Q}(-1-i^{2k},i^{2k})$ \[
\!\!\!\!\!\!\!\!\!\!\!\!\!\begin{array}{cc}
(u_{(i)})_{0}\leftrightarrow\left(\begin{array}{cc}
1 & \:0\\
0 & \:1\end{array}\right), & (u_{(i)})_{1}\leftrightarrow\left(\begin{array}{cc}
1 & \quad0\\
i^{2k-1}\sin\pi k^{} & \quad i^{2k}\end{array}\right),\\
(u_{(i)})_{2}\leftrightarrow\left(\begin{array}{cc}
0 & \:1\\
-i^{2k} & \:1+i^{2k}\end{array}\right), & (u_{(i)})_{3}\leftrightarrow\left(\begin{array}{cc}
\frac{1}{2}(1-i)(1+i^{2k}) & \quad i\\
i^{2k-1}\cos\pi k & \:\frac{1}{2}(1+i)(1+i^{2k})\end{array}\right).\end{array}\]
Reinterpreting, in accordance with (A.7) ((A.8)), $(u_{(1)})_{n}$
($(u_{(i)})_{n}$) as $((e_{(1)})_{n})^{\theta}$ (($(e_{(i)})_{n})^{\theta}$),
the following representations -- under the substitution $k=\theta=\tau$
-- hold true:\begin{equation}
\begin{array}{cc}
(e_{(1)})_{0}\leftrightarrow\left(\begin{array}{cc}
1 & \:0\\
0 & \:1\end{array}\right), & (e_{(1)})_{1}\leftrightarrow\left(\begin{array}{cc}
i & \quad0\\
2i\cos\frac{\pi}{2}\tau & \quad-i\end{array}\right),\end{array}\end{equation}
\[
\!\!\!\!\!\!\!\!\!\!\!\!\!\begin{array}{cc}
(e_{(1)})_{2}\leftrightarrow\left(\begin{array}{cc}
-\cot\frac{\pi}{2}\tau & \:\csc\frac{\pi}{2}\tau\\
-\csc\frac{\pi}{2}\tau & \:\cot\frac{\pi}{2}\tau\end{array}\right), & (e_{(1)})_{3}\leftrightarrow\left(\begin{array}{cc}
-i\cot\frac{\pi}{2}\tau & \quad i\csc\frac{\pi}{2}\tau\\
-i\csc\frac{\pi}{2}\tau+2i\sin^{}\frac{\pi}{2}\tau & \: i\cot\frac{\pi}{2}\tau\end{array}\right),\end{array}\]
\[
\begin{array}{cc}
(e_{(i)})_{0}\leftrightarrow\left(\begin{array}{cc}
1 & \:0\\
0 & \:1\end{array}\right), & (e_{(i)})_{1}\leftrightarrow\left(\begin{array}{cc}
1 & \quad0\\
2i^{\tau}\cos\frac{\pi}{2}\tau & \quad-1\end{array}\right),\end{array}\]
\begin{equation}
\begin{array}{cc}
(e_{(i)})_{2}\leftrightarrow\left(\begin{array}{cc}
-i\cot\frac{\pi}{2}\tau & \:1+i\cot\frac{\pi}{2}\tau\\
1-i\cot\frac{\pi}{2}\tau & \: i\cot\frac{\pi}{2}\tau\end{array}\right),\end{array}\end{equation}
\[
\begin{array}{cc}
(e_{(i)})_{3}\leftrightarrow\left(\begin{array}{cc}
-\cot\frac{\pi}{2}\tau & \quad-i+\cot\frac{\pi}{2}\tau\\
(-1+2\sin^{2}\frac{\pi}{2}\tau)(i+\cot\frac{\pi}{2}\tau) & \:\cot\frac{\pi}{2}\tau\end{array}\right).\end{array}\]
Up to poles for $\tau$ even, $(e_{(1)})_{n}$ and $(e_{(i)})_{n}$
all afford the required orthogonality and thus lay the basis for a
dual description of $\mathbf{Q}(0,1)$ and $\mathbf{Q}(0,-1)$. In
fact, the units $e'_{n}\equiv(e_{(1)})_{n}\vee i(e_{(i)})_{n}$ satisfy
the Hamilton relations for quaternions,\begin{equation}
\qquad(e')_{1}^{2}=(e')_{2}^{2}=(e')_{3}^{2}=e'_{1}e'_{2}e'_{3}=-e_{0},\end{equation}
while the units $\sigma_{n}\equiv(e_{(i)})_{n}\vee-i(e_{(1)})_{n}$
satisfy the Pauli relations for spin matrices,\begin{equation}
\qquad\begin{array}{c}
\sigma_{x}\sigma_{y}-\sigma_{y}\sigma_{x}=2i\sigma_{z},\\
\!\!\!\!\!\!\!\sigma_{x}\sigma_{y}+\sigma_{y}\sigma_{x}=0.\end{array}\qquad(x,y,z)=\mathrm{cycl}.\,(1,2,3)\end{equation}
\[
\quad\qquad\qquad\qquad\qquad\qquad\qquad\qquad\qquad\qquad\qquad\qquad\qquad\qquad\square\]

\end{document}